\documentclass[letterpaper, 10 pt, conference]{ieeeconf}

\IEEEoverridecommandlockouts                            
\usepackage{graphics} 
\usepackage{epsfig} 
\usepackage{times} 
\usepackage{amsmath} 
\usepackage{amssymb}  
\newtheorem{defi}{\textbf{Definition}}
\newtheorem{remark}{Remark}

\newtheorem{assume}{\textbf{Assumption}}
\newtheorem{lemma}{Lemma}
\usepackage[ruled,vlined]{algorithm2e}
\usepackage[utf8]{inputenc}
\usepackage{cite}
\usepackage{amsmath,amssymb,amsfonts}
\usepackage{algorithmic}
\usepackage{graphicx}
\usepackage{textcomp}
\usepackage{xcolor}
\usepackage{comment}
\usepackage{tikz}
\usepackage{accents}
\usepackage{multirow}

\newcommand*{\dt}[1]{%
   \accentset{\mbox{\large\bfseries .}}{#1}}
\newcommand*{\ddt}[1]{%
   \accentset{\mbox{\large\bfseries .\hspace{-0.25ex}.}}{#1}}
\def\BibTeX{{\rm B\kern-.05em{\sc i\kern-.025em b}\kern-.08em
    T\kern-.1667em\lower.7ex\hbox{E}\kern-.125emX}}

\title{\LARGE \bf
 A New Perspective of Accelerated Gradient Methods: The Controlled Invariant Manifold Approach}

\author{Revati Gunjal, Sushama Wagh, Syed Shadab Nayyer, 
Alex Stankovic, and Navdeep M. Singh
\thanks{Revati Gunjal, Sushama Wagh, Syed Shadab Nayyer, and Navdeep M. Singh are with Control and Decision Research Centre (CDRC), EED,
        Veermata Jijabai Technological Institute, Mumbai 400019, India.
        Sushama Wagh is a visiting scientist and Alex Stankovic is a distinguished scientist at SLAC National Accelerator Laboratory, Menlo Park, CA, USA.
         {\tt\small rgunjal\_p21@ee.vjti.ac.in}}%
 }

\begin{document}

\maketitle

\begin{abstract}
Gradient Descent (GD) is a ubiquitous algorithm for finding the optimal solution to an optimization problem. For reduced computational complexity, the optimal solution $\mathrm{x^*}$ of the optimization problem must be attained in a minimum number of iterations. For this objective, the paper proposes a genesis of an accelerated gradient algorithm through the controlled dynamical system perspective. The objective of optimally reaching the optimal solution $\mathrm{x^*}$ where $\mathrm{\nabla f(x^*)=0}$ with a given initial condition $\mathrm{x(0)}$ is achieved through control.

\end{abstract}
\begin{keywords}
Controlled Dynamical system, Gradient Descent, Manifold stabilization, Nesterov's Accelerated Gradient method, Passivity, and Immersion Approach
\end {keywords}
\section{Introduction}
Optimization plays a central role in various areas such as machine learning, data analysis, game theory, and a wide range of scientific and engineering applications. The gradient descent (GD) method is a ubiquitous algorithm for finding the optimal solution to an optimization problem through an iterative scheme. In the literature, the GD algorithm is analyzed as a first-order dynamical system, which helps one to get a better insight into the pros and cons of the algorithm. The gradient-based strategies for optimization are prevalent due to benefits like low iteration cost and easier deployment on parallel and distributed processing architecture. However, these methods are often inadequate to accommodate the growing complexity and scalability owing to large-scale data analysis that demands faster convergence rates. These issues motivated the requirement for accelerated gradient methods that achieve convergence to the optimal solution with the minimum number of iterations. B. T. Polyak in 1964 and  Y. Nesterov in 1983 developed a class of accelerated gradient algorithms with improved convergence rates and lower complexity as compared to the GD method. The notion of acceleration has also been incorporated with many optimization schemes such as stochastic optimization \cite{hu2009accelerated} \cite{lan2012optimal}, non-convex optimization \cite{ghadimi2016accelerated}\cite{li2015accelerated}, composite optimization \cite{beck2009fast}, and conic programming \cite{lan2011primal}, to name a few. These ideas have also been integrated into non-Euclidean optimization and higher-order optimization algorithms.

The literature discusses various continuous-time dynamical system formulations for accelerated gradient techniques, such as the Heavy Ball (HB) method \cite{polyak1964some}, Nesterov's Accelerated Gradient (NAG) \cite{nesterov1983method} method, and their extensions. The above approaches analyze the GD algorithm through a second-order ordinary differential equation (ODE). The momentum-based optimization algorithms, namely, HB and its extension Heavy Ball with Friction (HBF), are formulated based on the analogies of mechanical systems. The first-order GD system mimics the movement of a drop of water as it slides on the profile defined by the function $\mathrm{f(x)}$. The response of first-order GD is slow, and due to the inappropriate selection of step size, the algorithm might get stuck at the local minima. To tackle these limitations, B. T. Polyak introduced the HB algorithm \cite{polyak1964some} that accelerates the GD scheme by adding the momentum term. However, the response of the HB scheme exhibits an oscillatory behavior. To improve the response of the HB method, a damping term is added, establishing another strategy known as the HBF approach \cite{attouch2000heavy}. The ODE of HBF includes a ``\textit{geometric damping term}'' that involves the Hessian of the objective function.

Analysis of gradient-based techniques often incorporates the limiting ODE perspective. However, these ODEs cannot differentiate between the two fundamental algorithms, the HB method and Nesetrov's Accelerated Gradient method for Strongly Convex function (NAG-SC). In view of this, Bin Shi et. al. \cite{shi2021understanding} attempted an alternative limiting process to obtain the high-resolution ODEs for these methods. As a result of the high-resolution ODEs a geometric correction term appears. This geometric correction term includes Hessian of the function that provides information about the local geometry and the curvature of the descent. The knowledge of local geometry allows the selection of larger step sizes and hence faster convergence without risking the algorithm's stability or loss of any information.

As evident from the above-mentioned methods, the accelerated response can be achieved with fewer iterations by interpreting GD in terms of a second-order ODE. However, these methods throw no light on the origins of these algorithms and they still remain a conceptual mystery.
Our work attempts to provide insight into the origins of accelerated gradient algorithms and the approach so taken may be helpful in the design of algorithms rather than just providing an ODE-based interpretation of the existing accelerated gradient-based methods.
Heuristically, one could think of the algorithms as simply the numerical integration of a continuous time dynamical system. And there is a rich history of interpreting optimization algorithms as continuous time dynamical systems \cite{helmke2012optimization} \cite{brockett1991dynamical}\cite{soto2012nonlinear}. Further, L. Lassard et. al. \cite{lessard2016analysis} incorporated the control theoretic perspective in the accelerated gradient methods for optimization through the notion of Integral Quadratic Constraints (IQC). This was utilized for the stability analysis of the optimization algorithms to establish the upper bounds on the convergence rates of the gradient-based algorithms. Furthermore, a continuous time variational, continuous-time approach is proposed in \cite{wibisono2016variational} for understanding NAG method and its generalizations.
Motivated by the above approaches, this paper attempts to present the controlled dynamical system perspective to optimization algorithms and also explains the origin of the NAG in the continuous time domain.

The optimization problem is visualized as a manifold stabilization problem.
Manifold stabilization is a methodology developed in linear/nonlinear control theory that aims to design a control input or feedback scheme that guarantees the convergence (attractivity) of the system dynamics to a certain invariant manifold. The manifold selection depends on the desired objectives and if one considers the equilibrium point as the invariant manifold then the above problem reduces to the equilibrium stabilization problem.

The manifold stabilization in this work is achieved through the Passivity and Immersion (P\&I) approach \cite{9971723} \cite{nayyer2022passivity_gp}. This approach integrates the concepts of immersion intended for immersing the system dynamics into the defined lower-order target dynamics and the notion of attractivity of the manifold derived from the passivity theory. P\&I is proposed for control and stabilization, and the concept is extended to manifold stabilization. Recently the approach has been applied to parameter estimation \cite{nayyer2022passivity} by visualizing the estimation problem as a manifold stabilization problem, and the accurate estimation of the unknown parameters is accomplished via a control perspective. P\&I is an efficient and structured formulation that is useful to solve any problem by visualizing it as a manifold stabilization problem. An algorithm can be devised systematically through P\&I to obtain a controlled dynamical system that fulfills a specific objective.
This paper adapts the above approach to the optimization framework.

The controlled invariant manifold is susceptible to numerical errors, which can be viewed as perturbations, and as a result of which the invariant manifold may not persist when subjected to such perturbations. A novel approach is taken in this paper thereby the off-the-manifold dynamics transversal to the manifold is added to ensure stronger contraction of the ambient space and hence the persistence of the invariant manifold under numerical perturbations. This aids in designing a numerically insensitive algorithm for accelerated gradient descent.

The paper is structured as follows: The detailed procedure for the proposed controlled invariant manifold approach for the accelerated gradient method is elaborated in Section \ref{proposed methodology}. Section \ref{example}
provides insights into the applicability of the proposed theory. Section \ref{discussions} provides a brief discussion about the potential of the proposed philosophy. The paper is concluded along with the future work in Section \ref{conclusion}. Mathematical formulations of existing methods and a detailed procedure of the P\&I approach are described in the Appendix.

\section{Controlled Invariant Manifold Approach for Accelerated Gradient Method } \label{proposed methodology}
\subsection{Problem Definition}
Consider an unconstrained optimization problem 
    \begin{equation}\label{optimization problem}
    \begin{split}
  &\mathrm{min}\; \mathrm{f}(\mathrm{x})\\
  &\mathrm{s.t.}\; \mathrm{x}\in \mathbb{R} 
    \end{split}
\end{equation}
\begin{assume}\label{A1}
 The objective function $\mathrm{f}$ is smooth $\mu$-strongly convex with $L$-Lipschitz gradients. 
\end{assume}

\begin{defi}
    (Strong Convexity)
    A twice differentiable function $\mathrm{f:\mathbb{R}^n\rightarrow \mathbb{R}}$ is $\mu$-strongly convex \cite{wensing2020beyond} with $\mu>0$ being the modulus of convexity, if its Hessian matrix $\nabla^2\mathrm{f(x)}$  satisfies the matrix inequality
    \begin{equation}
        \nabla^2\mathrm{f(x)\geq \mu I\hspace{2cm}\forall x\in\mathbb{R}^n}
    \end{equation}
\end{defi}
\begin{defi}
    (L-Lipschitz Continuity)
    A differentiable function $\mathrm{f(x)}$ is smooth iff it has a Lipschitz continuous gradient, i.e. iff $\exists L<\infty$ such that
    \begin{equation}
        \left \|\nabla\mathrm{f}(\mathrm{x}_1)-\nabla\mathrm{f}(\mathrm{x}_2)  \right \|_2\leq L \left \| \mathrm{x}_1-\mathrm{x}_2 \right \|_2\: \: \: \: \forall \mathrm{x\in\mathbb{R}^n}
    \end{equation}
    
\end{defi}
The gradient descent (GD) for the above optimization problem can be written as: 
\begin{equation}\label{GD}
    \dt{\mathrm{x}}=-\nabla \mathrm{f(x)},
\end{equation}


The paper is concerned with a controlled dynamical system perspective to find the optimal point $\mathrm{x}^*$ of $\mathrm{f(x)}$ via a Controlled Invariant Manifold approach.

\subsection{Optimization as a Manifold Stabilisation Problem}

The concept of manifold stabilization has been incorporated into the optimization problem. In the context of control theory, a manifold is defined as a lower-dimensional sub-manifold immersed in the entire state space of the system. Instead of evolving in the whole state space the trajectories of the system dynamics evolve on the lower-dimensional manifolds. As per the desired behavior of the system or the set of constraints to be followed by the system, the manifold can be selected as an equilibrium point or a particular region of the state space of the system. 
 \subsubsection{Concept of Immersion}

Since the focus is on the manifold stabilization problem and its application in providing the controlled dynamical system perspective to the gradient method via the P\&I approach, the basic idea is to obtain an implicit manifold by projecting the NLS system  (system under consideration) onto a system with pre-defined characteristics. For example, in the equilibrium stabilization problem, the ultimate goal is to achieve the desired equilibrium in adaptive control and stabilization or to obtain a zero equilibrium in the state estimator design for the state estimation error.
Since the P\&I approach utilizes the notion of immersion from Immersion and Invariance (I\&I), the concept of immersing the system dynamics into the defined lower-order target dynamics is explained as follows.

Consider a NLS system 
\begin{equation}
 \mathrm{\dot{x}=f(x)+g(x)u}   
\end{equation}

where $\mathrm{x\in \mathbb{R}^n}$ and $\mathrm{u\in \mathbb{R}^m }$, the aim is to find control law $\mathrm{u=\nu (x)}$ such that stabilizes the closed loop system at an equilibrium point $\mathrm{x^*\in\mathbb{R}^n}$. 
\begin{itemize}
    \item \textbf{Target System}:
    The above system can be projected onto a target dynamical system 
\begin{equation} \label{TD}
\dt{\eta}=\beta(\eta)
\end{equation}
 where $\eta\in\mathbb{R}^\mathrm{p}$ and $\mathrm{p<n}$, having a GAS equilibrium at $\eta^*\in\mathbb{R}^\mathrm{p}$. 
 Under a smooth mapping $\pi: \mathbb{R}^\mathrm{p}\rightarrow \mathbb{R}^\mathrm{n}$,

 \begin{equation}
     \mathrm{x^*}=\pi(\eta^*)
 \end{equation}
    \item \textbf{Manifold}:
    The immersion of the system under consideration into lower-order target dynamics (\ref{TD}) via mapping $\mathrm{x}=\pi(\eta)$ helps in the construction of the implicit manifold $mathcal{M}$,
\begin{equation}
    \mathcal{M}=\left \{ \mathrm{x\in \mathbb{R}^n|x=\pi(\eta),\eta\in\mathbb{R}^p} \right \}
\end{equation}
    \item \textbf{Invariance of the Manifold}:
    The manifold $\mathcal{M=\left \{ \mathrm{x\in\mathbb{R}^n|\pi(x)=0} \right \}}$, with $\pi(\mathrm{x})$ being smooth, is said to be invariant, if $\mathrm{\pi(x(0))=0\Rightarrow \pi(x(t))=0, \; \;  \forall t\geq 0}$. i.e. all trajectories $\mathrm{x(t)}$ that start on the manifold remain on the manifold and asymptotically converge to the point $\pi(0)$. (Refer to Definition A.10 of \cite{astolfi2008nonlinear}) For proof of manifold invariance refer to \textit{Lemma 1}. 
    \item \textbf{Attractivity of the Manifold}:
    Once the invariant manifold is obtained the convergence of the off-the-manifold dynamics to the invariant manifold i.e. attractivity of the manifold is ensured by linking the passivity theory with immersion, which is accomplished by the P\&I approach discussed in the Appendix. 
\end{itemize}

\begin{figure}[ht!]
    \centering
    \includegraphics[width=\linewidth]{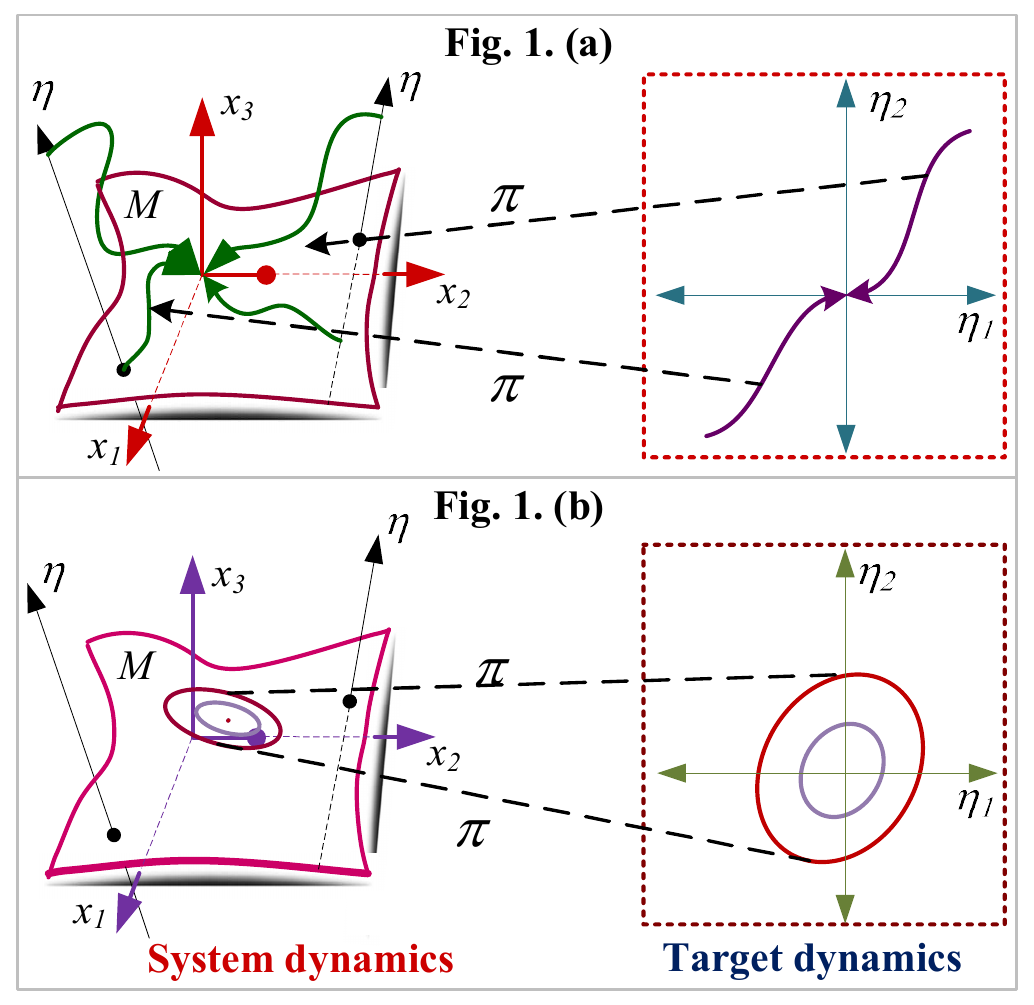}
    \caption{Schematic representation of the notion of immersion and its projection onto the target dynamics}
    \label{immersion 2}
\end{figure}

For further understanding of the immersion process and the choice of target dynamics according to the desired objective, a geometric visualization is depicted in Fig.\ref{immersion 2}. 
 Fig. \ref{immersion 2} illustrates two scenarios; equilibrium stabilization and orbital stabilization (with the goal is to induce or generate a periodic orbit in certain applications such as oscillation mechanisms in biology, electric motors,  power converters, and walking robots).

In Fig. \ref{immersion 2}(a), desired characteristics of the equilibrium stabilization for the third-order system (system under consideration) are captured by second-order target dynamics
  \begin{equation}\label{lOTG}
     \dt{\eta}_1=-\beta_1\eta_1 \hspace{1cm} \dt{\eta}_2=-\beta_2\eta_2.
 \end{equation}

Similarly,  Fig. \ref{immersion 2}(a) captures the desired characteristics of the periodic orbit for the third-order system by second-order target dynamics.
 \begin{equation}\label{odsoj}
     \dt{\eta}_1=\eta_2 \hspace{1cm} \dt{\eta}_2=-\mathrm{sin}(\eta_1)
 \end{equation}

\begin{remark}
    The objective of this article is to optimize the cost function. Thus, the manifold stabilization problem is viewed as an equilibrium stabilization problem, and focus is given to selecting the target dynamics  that have a GES equilibrium. 
\end{remark}

\subsubsection{Accelerated Gradient Method through P\&I}

Consider a second-order controlled dynamical system
\begin{equation} \label{control system}
           \dt{\mathrm{x}}_1=\mathrm{x_2} \hspace{1cm}
\dt{\mathrm{x}}_2=\mathrm{u}
   \end{equation}
   The control $\mathrm{u}\in \mathbb{R}$ is designed such that the trajectories of the above dynamical system will converge to the minima $\mathrm{x^*}$ of the objective function $\mathrm{f(x)}$. 
   
   $(S_1)$ \textbf{Construction of the implicit manifold:}

In the paper, the manifold is selected according to the desired objective, i.e.,  to achieve the optimal solution $\mathrm{x^*}$ of (\ref{optimization problem}) or the equilibrium point of the dynamical system in (\ref{control system}).  
The optimization problem in (\ref{optimization problem}) is geometrically visualized as the manifold stabilization problem to design an algorithm by adding the virtual control term $\mathrm{u}$. This input $\mathrm{u}$ is evaluated through the procedure of passivity and immersion (P\&I) approach \cite{nayyer2022towards} (The step-wise procedure is overviewed in Appendix) such that the trajectories converge exponentially to the equilibrium with a suitable choice of an implicit manifold.

  The target dynamics for finding the optimal solution $\mathrm{x^*}$ is defined as
 $\mathrm{x_2=-\beta\nabla f(x_1)}$. The target dynamics is defined on the implicit manifold 
   \begin{equation} \label{Implicit Manifold}
       \psi (\mathrm{x_1,x_2})=\mathrm{x_2+\beta\nabla f(x_1)}=0
   \end{equation}

\begin{lemma}
    The implicit manifold $\psi$ in (\ref{Implicit Manifold}) is invariant.
\end{lemma}
 \textit{Proof:} According to the relationship 
$\mathrm{x_2=-\beta\nabla f(x_1)}$, the manifold is obtained as $ \psi (\mathrm{x_1,x_2})=\mathrm{x_2+\beta\nabla f(x_1)}=0$. The normal to the velocity vector field $\begin{bmatrix}
\dot{\mathrm{x}}_1\\ 
\dot{\mathrm{x}}_2
\end{bmatrix}$ is given as
\begin{equation} \label{normal}
    \nabla \psi (\mathrm{x_1,x_2})=\begin{bmatrix}
 \beta\nabla^2\mathrm{f(x_1)} & \mathrm{I}
\end{bmatrix}
\end{equation}
The product of (\ref{normal}) and the velocity vector field yields $\beta\nabla^2 \mathrm{f(x_1)}\dot{\mathrm{x}}_1+\dot{\mathrm{x}}_2\Rightarrow \beta\nabla^2 \mathrm{f(x_1)}\dot{\mathrm{x}}_1-\beta\nabla^2 \mathrm{f(x_1)}\dot{\mathrm{x}}_1=0$. This implies that the velocity vector field is always tangent to the manifold $\psi$. Hence, the implicit manifold $\psi$ is invariant.

{$(S_2)$}  \textbf{Tangent space structure for control systems:}

  The psuedo Riemann metric $\mathrm{R(x)}=\mathrm{\nabla\psi (x_1,x_2)^T\nabla\psi (x_1,x_2)}$ is defined. 
  \begin{small}
   \begin{equation} \label{PR}
    \begin{split}
        \mathrm{R(x)} =\begin{bmatrix}
\mathrm{\beta \nabla^2f(x_1)}
\\ 
\mathrm{I}
\end{bmatrix} \begin{bmatrix}
\mathrm{\beta \nabla^2f(x_1)} & \mathrm{I}
\end{bmatrix}=\\
\begin{bmatrix}
 (\mathrm{\beta \nabla^2f(x_1)})^T (\mathrm{\beta \nabla^2f(x_1)})&\mathrm{\beta \nabla^2f(x_1)} \\ 
 \mathrm{\beta \nabla^2f(x_1)}& \mathrm{I}
\end{bmatrix}
    \end{split}
\end{equation}   
  \end{small}

Comparing (\ref{PR}) with (\ref{metric}), $\mathrm{m_{21}=\beta \nabla^2f(x_1)}$ and $\mathrm{m_{22}=1}$. Hence this gives $\mathrm{\nabla q(x)=\beta \nabla^2f(x_1)}$ and $\mathrm{q(x)=\beta \nabla f(x_1)}$. Since $\dt{\mathrm{x}}_2$ is along the vertical direction, the passive output is selected as $\dt{\mathrm{x}}_2+\mathrm{m}_{22}^{-1}\mathrm{m}_{21} \dt{\mathrm{x}}_1=\dt{\mathrm{x}}_2+\mathrm{\beta \nabla^2f(x_1)}  \dt{\mathrm{x}}_1$ which is in the same direction or parallel to $\dt{\mathrm{x}}_2$.

\textbf{$(S_3)$  Passive output:}
According to the preliminary definition of passive output, $\mathrm{y=y_1+y_2}$ is opted such that, $\mathrm{y_1}=\mathrm{\int_{0}^{t}(\nabla_{x_{1}}\varphi (x)\dt{\mathrm{x}}_1)dt}=\mathrm{\beta \nabla f(x_1)}$, and $\mathrm{y_2}=\mathrm{\int_{0}^{t}(\nabla_{x_{2}}\varphi (x)\dt{\mathrm{x}}_2)dt}=\mathrm{x_2}$.

$(S_4)$ \textbf{Storage function:}
The storage function (or candidate Lyapunov function) is defined with the passive output $\mathrm{y}$ as: 
\begin{equation} \label{storge func}
    \mathrm{S(x_1,x_2)=\frac{1}{2}(\mathrm{x_2+\beta \nabla f(x_1)})^2}
\end{equation}
The convergence of the trajectories of the off-the-manifold dynamics to the implicit manifold at an exponential rate $\alpha$ is ensured with the condition :
\begin{equation} \label{Expo conv}
    \mathrm{\dt{S}\leq -\alpha_1 S}
\end{equation}
Utilizing the condition (\ref{Expo conv}) along with the storage function (\ref{storge func}) yeilds 
\begin{small}
 \begin{equation}
    (\mathrm{x_2+\beta \nabla f(x_1)})(\dt{\mathrm{x}}_2+\mathrm{\beta \nabla^2f(x_1)}  \dt{\mathrm{x}}_1)\leq -\frac{\alpha_1}{2}(\mathrm{x_2+\beta \nabla f(x_1)})^2
\end{equation}   
\end{small}

 \begin{equation} \label{interm}
   \Rightarrow  \dt{\mathrm{x}}_2+\mathrm{\beta \nabla^2f(x_1)}  \dt{\mathrm{x}}_1+\frac{\alpha_1}{2}(\mathrm{x_2+\beta \nabla f(x_1)})=0
\end{equation}   

From (\ref{control system}), the equation (\ref{interm}) is modified in terms of the final control law as:

    \begin{equation}
       \mathrm{u}=  - \mathrm{\beta \nabla^2f(x_1)}\mathrm{x}_2-\frac{\alpha_1}{2}(\mathrm{x_2+\beta \nabla f(x_1)})
    \end{equation}

\begin{remark}
    For simplicity and understanding, the above formulation is illustrated for $\mathrm{x}\in \mathbb{R}$. This approach can be extended to $\mathrm{x}\in \mathbb{R}^\mathrm{n}$ with $\mathrm{u}\in \mathbb{R}^\mathrm{n}$.
\end{remark} 
\subsubsection{ODE of Accelerated Gradient Method through P\&I:}

Hence with the final control law, the controlled dynamical system modifies to the following form:
\begin{equation} \label{modified system}
\begin{split}
   \dt{\mathrm{x}}_1&=\mathrm{x_2}\\
\dt{\mathrm{x}}_2&=- \mathrm{\beta \nabla^2f(x_1)}\mathrm{x}_2-\frac{\alpha_1}{2}(\mathrm{x_2+\beta \nabla f(x_1)})
\end{split}
\end{equation}
The above system can be written as a second-order ordinary differential equation (ODE):
\begin{equation}\label{second order ODE}
    \ddt{\mathrm{x}}+\beta \nabla ^2 \mathrm{f(x)}\dot{\mathrm{x}}+\alpha \mathrm{\dot{x}}+\alpha\beta\nabla\mathrm{f(x)}=0
\end{equation}
where $\alpha=\frac{\alpha_1}{2}$.
\begin{remark}
    The ODEs of accelerated gradient methods in literature reflect a  geometric damping term involving the Hessian of the objective function. In P\&I the Hessian term $(\beta \nabla ^2 \mathrm{f(x)}\dot{\mathrm{x}})$ appears from the implicit manifold. The trajectories of the off-the-manifold dynamics take into account the geometry of the implicit manifold and approach the invariant manifold exponentially. The geometry of the manifold plays a crucial role in the acceleration phenomenon. Hence it justifies the term \textit{``clever geometric damping''} defined by \textit{Alvarez et. al.} \cite{alvarez2002second}. If the manifold is changed, the geometry will be altered, hence the Hessian damping term.  
\end{remark}

The above formulation is based on the invariant manifold theory which begins with the assumption that the invariant manifold exists. A further theory is developed from this point onwards by building upon this basic invariant manifold. From a numerical analysis viewpoint, it is necessary to analyze whether or not an invariant manifold persists \cite{talasila2006normally} under perturbation. And if so, if it maintains, loses, or gains differentiability is important.

For computational purposes, all the algorithms are implemented in the discrete-time system. These calculations are often subjected to numerical errors due to approximations. In order to enhance the accuracy in such circumstances the invariant manifold should persist under the perturbations (i.e., numerical errors).

\subsection{Persistence of Invariant Manifold}
If the considered invariant manifold controlled with P\&I possesses the property of normal hyperbolicity, then the normally hyperbolic invariant manifold (NHIM) persists under perturbations. The notion of normal hyperbolicity and the conditions to analyze its persistence property is introduced by N. Fenichel \cite{fenichel1971persistence} and M. Hirsch, et.al. \cite{hirsch1977invariant}.
 For a dynamical system, an invariant manifold is considered as NHIM \cite{capinski2020persistence} if there is a splitting of state space into three invariant sub-bundles: a bundle tangent to the manifold, a stable bundle (dynamics on this is contracting), and an unstable bundle (dynamics on this is expanding). If there exists a dominance of dynamics on the stable and unstable bundle over the dynamics on the tangent bundle, then the manifold is NHIM. 
\begin{remark}
     The normal hyperbolicity implies that the dynamics transverse to the invariant manifold is dominant i.e., the growth rate of the vectors transverse to the manifold dominates the growth rate of vectors tangent to the manifold. This conceptual explanation is understood from the definition of normal hyperbolicity and $r$-normally hyperbolicity (Refer \cite{hirsch1977invariant} )
\end{remark}
\subsubsection{Concept of Normally Controlled Hyperbolicity}
Since the concerned invariant manifold is a controlled invariant manifold hence the notion of $r$-normally controlled hyperbolicity introduced in \cite{talasila2006normally} should be considered.

Given a controlled invariant manifold $\mathbb{M}$, stabilised by the control law $\mathrm{u}$ designed for the controlled dynamical system $\mathrm{\dt{x}=X_C(x,u)}$, where $\mathrm{x\in \Omega }$ and $\mathrm{u}\in \mathbb{U}$. The corresponding flow for the controlled invariant manifold is $f:\mathbb{M}\times \mathbb{U}\rightarrow \mathbb{M}$. Assume $\sigma _\mathrm{C}$ chosen to associate to a state $\mathrm{x\in \Omega }$, a stabilising control law $\mathrm{u}\in \mathbb{U}$. 
\begin{defi}
    (r-normally controlled-hyperbolicity) \\
    In the above-mentioned controlled dynamical system, the manifold $\mathbb{M}$ is an r-normally controlled-hyperbolic with respect to the controlled flow
    \begin{equation*}
        \rho := f\circ \sigma _\mathrm{C}
    \end{equation*}
\end{defi}
iff the tangent bundle of $\Omega$ restricted to $\mathbb{M}$ splits into two continuous sub-bundles:
\begin{equation*}
    \mathbb{T}_\mathbb{M}\Omega =\mathbb{T}\mathbb{M}\oplus \mathbb{N}^{\mathrm{st}}
\end{equation*}
invariant by the tangent of $\rho$ i.e. $\mathbb{T}{\rho}$, such that for $p\in\mathbb{M}$, $0\leq \kappa \leq \mathrm{r}$:
\\ $\mathbb{T}\rho$ contracts the stable normal bundle $\mathbb{N}^{\mathrm{st}}$ more sharply than $\mathbb{T}\mathbb{M}$, i.e.,
\begin{equation*}
  \left \| (\mathbb{T}\rho\setminus \mathbb{N}_p^{\mathrm{un}}(\mathbb{M}))^{-1} \right \|\cdot \left \| \mathbb{T}\rho\setminus \mathbb{N}^{\mathrm{st}}(\mathbb{M}) \right \|^\kappa\leq c(1/\lambda )  
\end{equation*}
for some $1<\lambda <\infty, 0<c<\infty$.
Referring to the contraction and expansion conditions, 
the concept of normal hyperbolicity can be informally defined as the normal behavior of $\mathbb{T}\rho$ dominates the tangential behavior of $\mathbb{T}\rho$. The dominance of the normal behavior over the tangential behavior allows the controlled invariant manifold to persist subjected to the perturbations. 
\subsubsection{ ODE of Accelerated Gradient Method through P\&I and Additional off-the-manifold Dynamics}
To ensure the persistence of the invariant manifold under the perturbations, a contracting off-the-manifold dynamics transverse to the invariant manifold having the same equilibrium point as (\ref{second order ODE}) is added. Hence the equation (\ref{second order ODE}) modifies to the following form:  
\begin{equation}\label{PNI+Forcing}
   \ddot{\mathrm{x}}+\beta \nabla ^2 \mathrm{f(x)}\dot{\mathrm{x}}+\alpha \mathrm{\dot{x}}+\alpha\beta\nabla\mathrm{f(x)}+\left [ \alpha\dot{\mathrm{x}}+\nabla\mathrm{f(x)} \right ]=0
\end{equation}

The dynamics $\mathrm{\dt{x}=-\frac{1}{\alpha}\nabla f(x)}$ transverse to the invariant manifold aids to the effect of contraction of the stable normal bundle $\mathbb{N}^{\mathrm{st}}$. In effect, the dynamical behavior normal to the invariant manifold is more dominant than the form of the dynamics defined on the invariant manifold. 
The additional off-the-manifold dynamics adhere to the notion of normal hyperbolicity of the controlled invariant manifold. i.e. the controlled invariance properties are preserved under perturbations. 
\begin{remark}
  In the considered problem, the optimal solution $\mathrm{x^*}$ is hyperbolic (Refer \cite{wang2015immersion}). Under normal hyperbolicity, even if the system perturbs, the new point of a perturbed system preserves the properties of the interest. 
\end{remark}

Due to the inclusion of the dynamics transverse to the invariant manifold, the space around the manifold becomes more contracting. The strong contraction ensures the persistence of the invariant manifold and results in faster convergence of the trajectories of the off-the-manifold dynamics towards the invariant manifold. The persistence of the invariant manifold under the perturbation ensures that the proposed algorithm is insensitive to numerical errors. 
\begin{remark}
The stability of the system remains unaffected by the addition of the dynamics transverse to the manifold since the response of the additional dynamics is exponential and its effect diminishes as the trajectories of the off-the-manifold dynamics approach the equilibrium. Since the effect of dynamics transverse to the manifold disappears at the equilibrium point, there is no shift of the equilibrium point.
\end{remark}



\subsection{Correlation with the Nesterov's accelerated gradient method }
The ODE derived with manifold stabilization procedure
\begin{equation}\label{Proposed ODE}
    \ddt{\mathrm{x}}+2\alpha \mathrm{\dt{x}}+\beta \nabla ^2 \mathrm{f(x)}\dt{\mathrm{x}}+(1+\alpha\beta)\nabla\mathrm{f(x)}=0
\end{equation}
is similar to Nesterov's Accelerated Gradient for strongly convex functions (NAG-SC) mentioned in the high-resolution differential equation derived by \textit{Bin Shi et. al.}  \cite{shi2021understanding}).
\begin{equation}\label{NAG-sc}
    \ddt{\mathrm{x}}+2\sqrt\mu \mathrm{\dt{x}}+\sqrt s \nabla ^2 \mathrm{f(x)}\dt{\mathrm{x}}+(1+\sqrt {\mu s})\nabla\mathrm{f(x)}=0
\end{equation}

The Hessian of the function in the above equation works as the geometric damping and aids in the acceleration of the GD algorithm. The algorithms designed in the continuous time system are implemented with the help of numerical integration in the discrete-time system for computation purposes. To reduce the computational burden on the processor, the computation complexity should be reduced. To achieve reduced complexity, the selection of the step size is a  matter of concern.
The selection of the larger step size allows the reduction in the number of iterations. But the larger step size poses a serious problem to the stability of the algorithm.  The Hessian of the function gives information about the curvature of the descent resulting in faster convergence by appropriate scaling of the step size as per the level of the curvature.

The Hessian term reflects in the P\&I approach by following a systematic methodology without any approximations. The Hessian damping term in the ODE derived through the P\&I is governed by the geometry of the manifold and allows the convergence to the solution with a larger step size without disturbing the stability of the algorithm. Larger step sizes ensure the convergence to the optimal solution with a lesser number of iterations leading to accelerated algorithms.

\subsection{Significance of $\alpha$ and $\beta$}
The parameters $\alpha$ and $\beta$ in the proposed approach can be selected by comparing the ODE (\ref{Proposed ODE}) to the parameters $\mu$ and $s$ in the NAG-sc. The parameters $\mu$ (the modulus of convexity of the strongly convex function) and $s$ (step size) appear from the discretization. These parameters play an important role in the stability of the algorithm and are selected by following several guidelines. 

The parameter $\alpha$ present in the off-the-manifold dynamics is bounded by the modulus of convexity of the function. 
$\alpha$ is related to the manifold attractivity (the trajectories of the off-the-manifold dynamics approach the implicit manifold exponentially with the rate $\alpha$) and $\mu$ represents the lower bound to the growth of the function. The trajectories of the off-the-manifold dynamics should approach the invariant manifold faster than this bound. Hence $\alpha$ can be related to $\mu$. 
The other parameter $\beta$ present in the dynamics defined on the invariant manifold is bounded by the Lipschitz constant of the function. $\beta$ in the definition of the implicit manifold is similar to the step size $\mathrm{s}$ in the GD iterative algorithm. The optimal step size is bounded by the Lipschitz constant of the function (Refer to equation (1.2.12) in \cite{nesterov2003introductory}). 
Referring to the NAG-SC in equation (\ref{NAG-sc}), the optimal values of the parameters are selected as $\alpha=\sqrt{\mu}$ and $\beta=\sqrt{s}$. 
\subsection{Stability Analysis}
The ODE (\ref{PNI+Forcing}) for the accelerated algorithm can be looked upon as a second-order linear dynamical system as follows:
\begin{equation}
    \begin{bmatrix}
\dt{\mathrm{x}}_1\\ 
\dt{\mathrm{x}}_2
\end{bmatrix}=\begin{bmatrix}
0 &1 \\ 
 0&0 
\end{bmatrix}\begin{bmatrix}
\mathrm{x_1}\\ 
\mathrm{x_2}
\end{bmatrix}+\begin{bmatrix}
0 &0 \\ 
1 & 1
\end{bmatrix}\begin{bmatrix}
\mathrm{u_1}\\ 
\mathrm{u_2}
\end{bmatrix}
\end{equation}
Since the dynamical system is linear, by superposition property, both the control inputs in the above system $\mathrm{u_1}$ and $\mathrm{u_2}$ can be evaluated independently. Control input $\mathrm{u_1}$ is evaluated through the P\&I procedure, whereas the other control input $\mathrm{u_2}$ is added for the persistence of the invariant manifold under perturbations. Hence the additional off-the-manifold dynamics transversal to the invariant manifold will formulate the control input $\mathrm{u_2}=\mathrm{-\alpha x_2-\nabla f(x_1)}$.

The Lyapnouv function for the system dynamics with $\mathrm{u_1}$ evaluated through the P\&I strategy is mentioned earlier in equation (\ref{storge func}). The convergence of the trajectories of the off-the-manifold dynamics to the invariant manifold with the control input $\mathrm{u_1}$ is ensured with the condition (\ref{Expo conv}).

The stability of the system with control input $\mathrm{u_2}$ can be evaluated with $\mathrm{u_1}=0$. The second order dynamical system with control input $\mathrm{u_2}$ will be of a form:
\begin{equation} \label{second input}
\begin{split}
   \dt{\mathrm{x}}_1&=\mathrm{x_2}\\
\dt{\mathrm{x}}_2&=- \mathrm{\alpha x_2- \nabla f(x_1)}
\end{split}
\end{equation}

For the sake of simplicity and understanding, the stability analysis is illustrated for $\mathrm{x_1,x_2\in\mathbb{R}}$.
A similar analysis can be suitably extended to $\mathrm{x_1,x_2\in\mathbb{R}^n}$.
\begin{lemma}
   The second-order dynamics in (\ref{second input}) is stable.
\end{lemma}
 \textit{Proof:} The Lyapunov function for the second order system mentioned in (\ref{second input}) with the additional off-the-manifold dynamics transversal to the invariant manifold is given as:
 \begin{equation}
     \mathrm{V(x_1,x_2)}=\mathrm{f(x_1)-f(x_1^*)+\frac{1}{2}{x}_2^2}
 \end{equation}
 For proving stability, a derivative of the above Lyapunov function is evaluated as:
 \begin{equation}
     \mathrm{\dt{V}(x_1,x_2)}=\mathrm{\nabla f(x_1)\dt{x}_1+x_2\dt{x}_2}
 \end{equation}
 \begin{equation}
     \Rightarrow \mathrm{\dt{V}(x_1,x_2)}=\mathrm{\nabla f(x_1)x_2-\alpha x_2^2-\nabla f(x_1)x_2}
 \end{equation}
 \begin{equation}
     \Rightarrow \mathrm{\dt{V}(x_1,x_2)}=-\alpha \mathrm{x_2^2}
 \end{equation}
 Since $\mathrm{\dot{V}(x_1,x_2)}\leq 0$ the second-order system dynamics (\ref{second input}) is stable.

\begin{remark}
    The asymptotic stability of the above system can be proved using LaSalle's theorem. 
\end{remark}

\begin{remark}
    For the exponential stability of the system in (\ref{second input}), the Lyapunov function from \cite{shi2021understanding} can be referred. From the Theorem 2 of \cite{shi2021understanding} the Lyapunov function can be selected as 
    \begin{equation*}
       \mathrm{V(x_1,x_2)}=\mathrm{f(x_1)-f(x_1^*)+\frac{1}{2}{x}_2^2+\frac{1}{2}(x_2+\alpha(x_1-x_1^*))^2} 
    \end{equation*}
    where $\alpha=\alpha_1 /2$.

    By following the same steps as mentioned in \cite{shi2021understanding}, the exponential stability can be proved with the help of the above Lyapunov function. 
    The exponential convergence rates for both the dynamics with $\mathrm{u_1}$ and $\mathrm{u_2}$ are the same i.e. $\alpha_1$. 
\end{remark}

\section{Numerical Example and Results}
\label{example}
To gain deeper insights into the proposed approach consider an optimization problem with a quadratic cost function.
    \begin{equation}
    \begin{split}
  &\mathrm{min}\; \mathrm{f}(\mathrm{x})\\
  &\mathrm{s.t.}\; \mathrm{x}\in \mathrm{X}
    \end{split}
\end{equation}
where   $\mathrm{X}=\left \{ \mathrm{x}\in \mathbb{R} \right \}$ and $\mathrm{f(x)=\frac{1}{2}x^Tx}$.
   
    The gradient descent for the above quadratic objective function will be:
    \begin{equation}
        \dt{\mathrm{x}}=-\mathrm{\nabla f}=-\mathrm{x}
    \end{equation}
    For solving the given problem with the control-theoretic approach consider a second-order system
   \begin{align}
\begin{split}
\mathrm{\dt{x}}_1=\mathrm{x_2}\hspace{1cm}
\mathrm{\dt{x}}_2=\mathrm{u} 
\end{split} 
\end{align}
As per the equation (\ref{Implicit Manifold}), the implicit manifold is defined as 
\begin{equation}
     \psi (\mathrm{x_1,x_2})=\mathrm{x_2+\beta x_1}=0
\end{equation}
Following the remaining steps of P\&I the control law is derived.
\begin{equation}
    \mathrm{u}=\mathrm{-\beta x_2-\alpha x_2-\alpha\beta x_1}
\end{equation}
The corresponding second-order ODE is
\begin{equation} \label{ODE_pni}
    \mathrm{\ddt{x}}+\beta \mathrm{\dt{x}}+\alpha \mathrm{\dt{x}}+\alpha \beta \mathrm{x}=0
\end{equation}
\subsection{Response before the addition of dynamics transversal to the invariant manifold}
\begin{figure}[ht] \label{pni_alpha_10}
 \centering
     \includegraphics[ width=\linewidth]{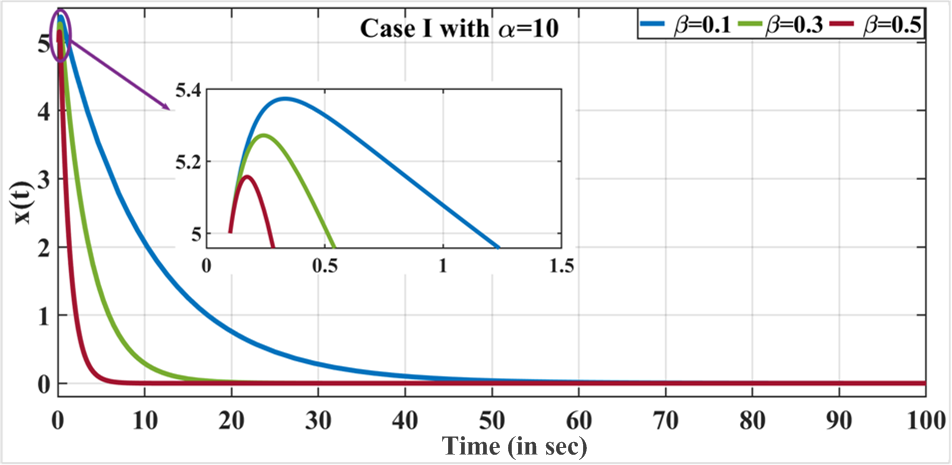}
     \caption{Response of Accelerated Gradient Descent derived through P\&I with the y-axis representing the solution of the quadratic cost function $\mathrm{f(x)}$ and x-axis representing the time for $\alpha=10$ with various values of $\beta$  **}
     
 \end{figure}

 \begin{figure}[ht] \label{pni_alpha_1}
 \centering
     \includegraphics[ width=\linewidth]{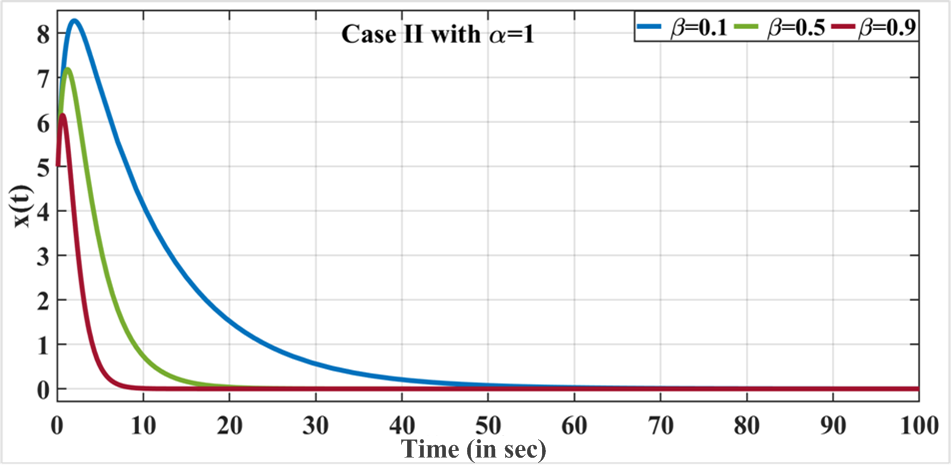}
     \caption{Response of Accelerated Gradient Descent derived through P\&I with the y-axis representing the solution of the quadratic cost function $\mathrm{f(x)}$ and x-axis representing the time for $\alpha=1$ with various values of $\beta$  **}  
 \end{figure}
 The ODE (\ref{ODE_pni}) is simulated with the various values of the parameters $\alpha$ and $\beta$ as demonstrated in Fig. 2 and Fig. 3. It is observed that the convergence to the optimal solution $\mathrm{x^*}$ becomes faster with the increasing value of the parameter $\beta$. From Fig. 2, it is visible that the peak overshoots are decreased as the value of $\alpha$ is increased as compared to Fig. 3. 

\subsection{Response after the addition of dynamics transversal to the invariant manifold}

For ensuring the persistence of the invariant manifold under the perturbations such as numerical errors, another dynamics transverse to the invariant manifold is added. Hence the ODE (\ref{ODE_pni}) modifies as:
\begin{equation} \label{ODE_overall}
    \mathrm{\ddt{x}}+\beta \mathrm{\dt{x}}+2\alpha \mathrm{\dt{x}}+(1+\alpha \beta) \mathrm{x}=0
\end{equation}
  \begin{figure}[ht]
 \centering
     \includegraphics[ width=\linewidth]{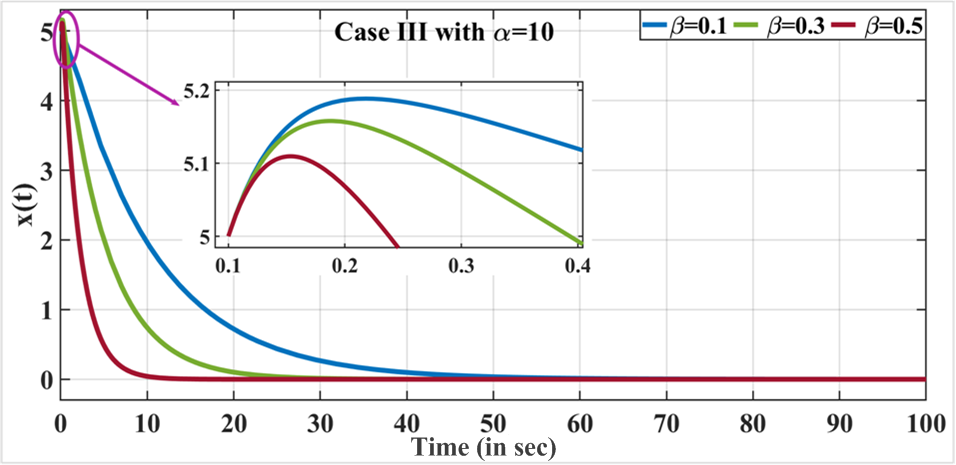}
     \caption{Effect of the addition of another dynamics transverse to the invariant manifold to the Accelerated Gradient Descent derived through P\&I with the y-axis representing the solution of the quadratic cost function $\mathrm{f(x)}$ and x-axis representing the time for $\alpha=10$ with various values of $\beta$  **}
     \label{AGD_alpha_10}
 \end{figure}
   \begin{figure}[ht]
 \centering
     \includegraphics[ width=\linewidth]{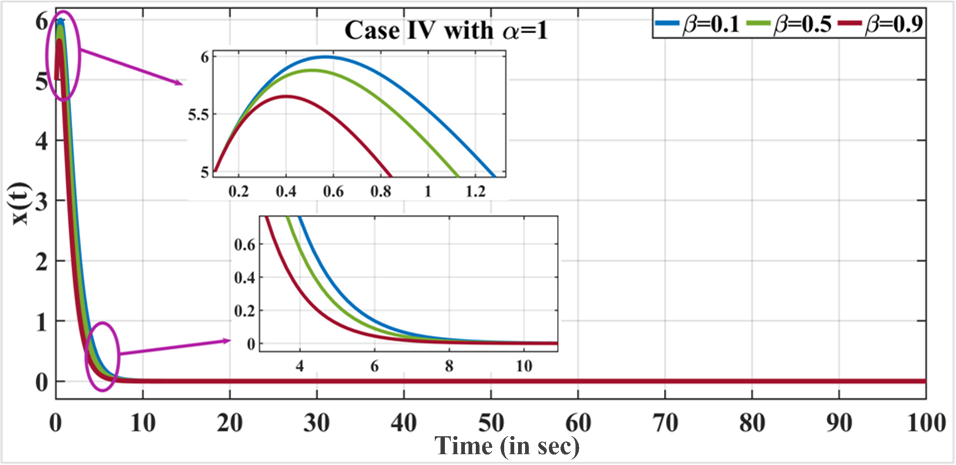}
     \caption{Effect of the addition of another dynamics transverse to the invariant manifold to the Accelerated Gradient Descent derived through P\&I with the y-axis representing the solution of the quadratic cost function $\mathrm{f(x)}$ and x-axis representing the time for $\alpha=1$ with various values of $\beta$  **}
     \label{AGD_alpha_1}
 \end{figure}

The above ODE is simulated for the different values of $\alpha$ and $\beta$ similar to the earlier case and the results are illustrated in Fig. \ref{AGD_alpha_10} and Fig. \ref{AGD_alpha_1}. For the given quadratic cost function the optimal values of the parameters are selected as $\alpha=1$ and $\beta=0.9$. 

The effect of the additional off-the-manifold dynamics transverse to the invariant manifold with the optimal values of parameters $\alpha$ and $\beta$ is demonstrated in Fig. \ref{AGD_alpha_1}. It is observed that due to the dominating effect of the dynamics transverse to the invariant manifold over the dynamics defined on the invariant manifold, the convergence to the optimal solution $\mathrm{x^*}$ of the quadratic cost function is achieved with the minimum number of iterations. 
\begin{remark}
   **  Fig. 2 to Fig. 5 are revised in the latest version of the manuscript. The ODE response is illustrated in continuous time domain and the axes labels are revised accordingly.
    \end{remark}
\section{Discussions} \label{discussions}
In this paper, a controlled dynamical system perspective for deriving the accelerated gradient method is proposed. The perspective is motivated by an understanding that the algorithms are numerical integration of a dynamical system. Hence for achieving a particular objective, any problem could be looked upon as a dynamical system. And a control law is designed for the dynamical system to achieve the desired objective. The controlled dynamical system perspective is implemented with the incorporation of the notion of manifold stabilization. A control scheme for the dynamical system is designed to ensure the convergence (attractivity) of the system dynamics to a particular invariant manifold. Hence the considered problem is also referred to as the controlled invariant manifold approach. Since the objective of this paper is to optimally reach the optimal solution $\mathrm{x^*}$ where $\mathrm{\nabla f(x)=0}$, the manifold is considered as an equilibrium point. Hence the considered problem is reduced to an equilibrium stabilization problem.

The equilibrium stabilization is achieved by blending the concept of immersion and passivity through the P\&I approach. The invariant manifold is constructed through the process of immersion by immersing the system dynamics into the specified lower-order target dynamics. The convergence of the trajectories of the system dynamics to the invariant manifold (attractivity of the invariant manifold) is derived through the passivity theory. The ODE derived through the systematic procedure of P\&I contains a term with Hessian of the function which provides the knowledge of the local geometry of the curvature and enables the convergence to the optimal solution in fewer iterations with the selection of larger step sizes.  

Since the algorithms are implemented through the various techniques of numerical integration in the discrete-time system, some numerical errors occur due to the approximations involved in the process. The invariant manifold is often susceptible to these numerical errors (perturbations) and the invariant manifold might not persist under such perturbations. The proposed theory is developed on the basic invariant manifold and this invariant manifold should persist when subjected to numerical errors. Hence, to ensure the persistence of the controlled invariant manifold, an additional contracting dynamics transversal to the manifold is added which constitutes another major contribution of this paper. The additional off-the-manifold dynamics adhere to the notion of normal hyperbolicity of the controlled invariant manifold. Under normal hyperbolicity, even if the system perturbs, the new point of the perturbed system preserves the properties of the interest. The additional dynamics ensure the stronger contraction of the ambient space and hence the persistence of the invariant manifold under numerical perturbations. As a result, a numerically insensitive algorithm for the accelerated gradient method is developed. 

It is observed that the ODE obtained through the proposed approach with the additional off-the-manifold dynamics resembles the high-resolution ODE for the NAG-SC method. The basic formulation derived for the NAG-SC could be modified to other settings for accelerated gradient methods. For various objective functions, the Hessian of the function term is negligible and could be omitted. This leads to the equation which is similar to the high-resolution ODE for the HB algorithm. Since the proposed formulation is developed on the persistent invariant manifold, the algorithm remains insensitive to the additional terms to the dynamics. Hence with the addition of another damping term, the proposed formulation could be extended to the accelerated Triple Momentum (TM) method \cite{sun2020high}. These extensions are summarised with the corresponding ODEs in the following figure. 
\begin{figure}[ht] \label{equations}
 \centering
     \includegraphics[ width=\linewidth]{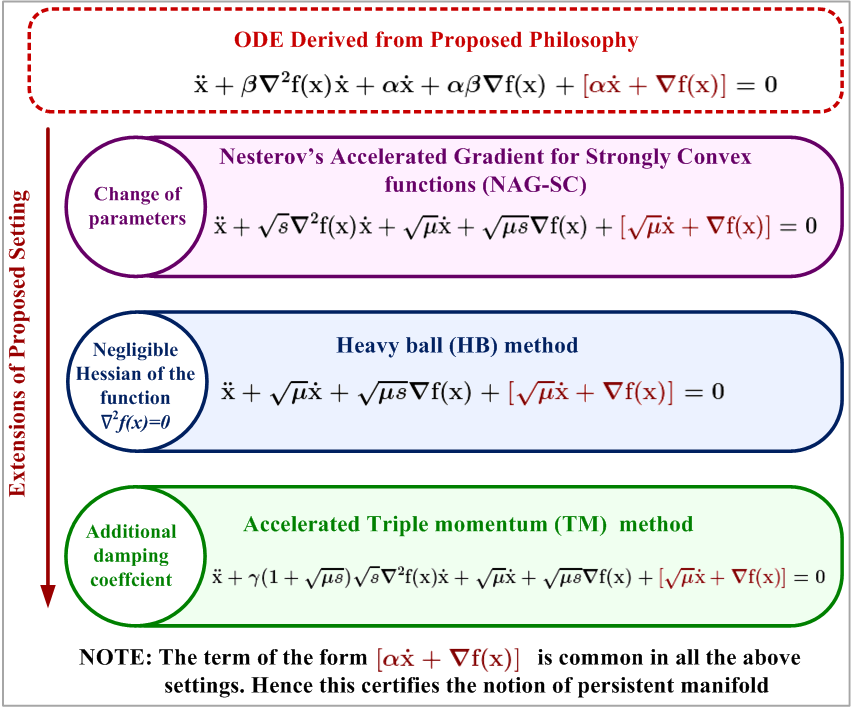}
     \caption{Extensions of ODE obtained through proposed approach }
     
 \end{figure}
 The term $[\alpha \dt{\mathrm{x}}+\nabla \mathrm{f(x)}]$ appears in NAG-SC, HB and TM methods. This certifies that the additional dynamics is necessary for the accelerated gradient methods. By combining the additional dynamics with the remaining ODEs the above equations are rewritten as:
  \begin{itemize}
     \item \textbf{Proposed ODE:}
      \end{itemize}
    \begin{equation*}
         \ddt{\mathrm{x}}+2\alpha \mathrm{\dt{x}}+\beta \nabla ^2 \mathrm{f(x)}\dt{\mathrm{x}}+(1+\alpha \beta)\nabla\mathrm{f(x)}=0
     \end{equation*}
 \begin{itemize}
     \item \textbf{Nesterov's Accelerated Gradient for Strongly Convex function (NAG-SC):}
      \end{itemize}
     \begin{equation*}
         \ddt{\mathrm{x}}+2\sqrt\mu \mathrm{\dt{x}}+\sqrt s \nabla ^2 \mathrm{f(x)}\dt{\mathrm{x}}+(1+\sqrt {\mu s})\nabla\mathrm{f(x)}=0,
     \end{equation*}
     where the parameters of the proposed ODE are selected as $\alpha=\sqrt{\mu}$ and $\beta=\sqrt{s}$.
\begin{itemize}
     \item \textbf{Heavy ball method (HB):}
      \end{itemize}
     \begin{equation*}
         \ddt{\mathrm{x}}+2\sqrt\mu \mathrm{\dt{x}}+(1+\sqrt {\mu s})\nabla\mathrm{f(x)}=0
     \end{equation*}
        where the parameters of the proposed ODE are selected as $\alpha=\sqrt{\mu}$ and $\beta=\sqrt{s}$ and $\nabla ^2 \mathrm{f(x)}=0$.
\begin{itemize}

     \item \textbf{Accelerated Tripple Momentum method (TM):}
      \end{itemize}
     \begin{equation*}
        \ddt{\mathrm{x}}+2\sqrt\mu \mathrm{\dt{x}}+\gamma(1+\sqrt {\mu s})\sqrt s \nabla ^2 \mathrm{f(x)}\dt{\mathrm{x}}+(1+\sqrt {\mu s})\nabla\mathrm{f(x)}=0
     \end{equation*}
        where the parameters of the proposed ODE are selected as $\alpha=\sqrt{\mu}$ and $\beta=\sqrt{s}$ and another damping coefficient $\gamma$ is added. 
     \begin{remark}
   The stability certificate of these methods can be referred from the proofs given in \cite{shi2021understanding} \cite{sun2020high}.       
     \end{remark}
\section{Conclusion and Future Work}
 The controlled invariant manifold approach aims to promote a novel philosophy to provide a better understanding of the working of the algorithm. The proposition forwarded by this work is problems in optimization could be solved by visualizing them as a manifold stabilization problem. And the algorithm could be generated systematically through the P\&I approach to obtain a controlled dynamical system that achieves a desired objective. The authors are intending to extend the proposed approach to the geodesic convexity settings. In future work, the above setting proposed for the unconstrained optimization case will be extended to the primal-dual optimization problems. 

\label{conclusion}
\section*{Acknowledgment}
The authors of the paper would like to acknowledge the support provided by Control and Decision Research Center (CDRC), VJTI, Mumbai, India, and (Project 20627-X2001, FY22-Stankovic Fusion Methods) SLAC National Accelerator Laboratory, Menlo Park, CA, USA.

\section*{Appendix}
\subsection{Existing Approaches for Accelerated Gradient Method}\label{Appendix 1}

Consider an unconstrained optimization problem of the form 
\begin{equation*}
    \mathrm{min\;  f(x)} \;\;\;\;\;\; \mathrm{x\in \mathbb{R}^n}
\end{equation*}

The objective function $\mathrm{f(x)}$ can be a convex function or a strongly convex function. The optimal value of this function is denoted as $\mathrm{f^*=\underset{x}{min}\; f(x)}$ and the optimal solution of the above-mentioned problem is represented as $\mathrm{x^*=arg\:\underset{x}{ min}\; f(x)}$. The straightforward approach for finding $\mathrm{x^*}$ is to solve the system of equation $\mathrm{\nabla f(x^*)=0}$ which is usually not easier to solve from an analytical viewpoint. Hence, a preferable strategy is an iterative scheme that involves computing a minimizing sequence of points $\mathrm{x^{(0)}, x^{(1)}, \cdots }$ such that $\mathrm{f(x^{(k)})\rightarrow f(x^{*})  }$ as $\mathrm{k\rightarrow \infty }$.

Gradient descent (GD) or steepest descent is a ubiquitous iterative algorithm for solving optimization problems. It is represented  by the first-order dynamics in the time system as follows:
\begin{equation}
    \begin{split}
        \mathrm{\dt{x}(t)+\nabla f(x(t))=0}\\
\mathrm{x(0)=x_0}
    \end{split}
\end{equation}
Since GD is a first-order system, its trajectory is solely governed by the initial point $\mathrm{x_0}$.
The GD system mimics the movement of a drop of water as it slides on the profile defined by function $\mathrm{f(x)}$. The GD algorithm is slower and it often gets stuck in local minima due to inappropriate selection of step size.  

The limitations of GD are overcome by the Heavy Ball (HB) method \cite{polyak1964some} proposed by B. T. Polyak which is the accelerated version of the GD represented as:
\begin{equation}
    \mathrm{\ddt{x}(t)+\nabla f(x(t))=0}
\end{equation}
HB method provides an accelerated response with oscillations. To suppress these oscillations H. Attouch et. al. \cite{attouch2000heavy} added a damping term  representing friction and the modified method is referred to as Heavy Ball with Friction (HBF) system. It is a non-linear oscillator with damping, given as:
\begin{equation}
    \mathrm{\ddt{x}(t)+\lambda \dt{x}(t)+\nabla f(x(t))=0}
\end{equation}

\subsection{The P\&I Approach}\label{Appendix 2}
The immersion of the system under consideration into lower-order target dynamics helps in the construction of the implicit manifold.
The convergence and  manifold attractivity are ensured by linking the passivity theory with immersion. The P\&I approach utilizes the idea of an invariant target manifold giving rise to a non-degenerate two form, through which the definition of certain passive outputs and storage functions leads to a generation of control law for stabilizing the system. Moreover, the involvement of a particular structure of tangent bundle associated with a suitable pseudo-Riemannian structure imposed on the control system makes the P\&I more methodical and systematic.

\begin{figure}[ht!]
    \centering
    \includegraphics[width=\linewidth]{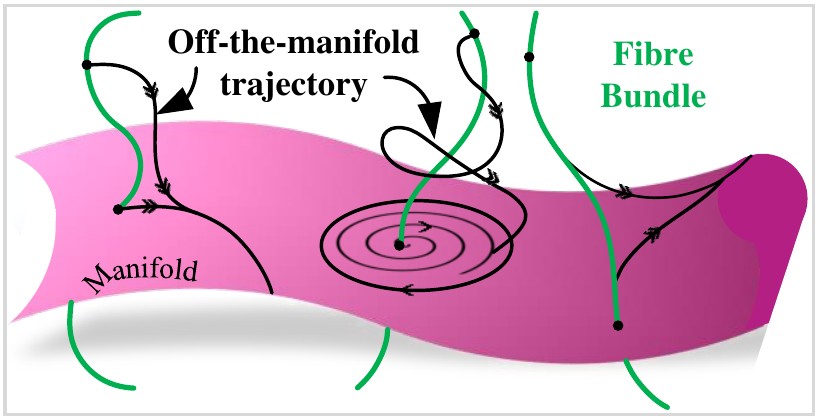}
    \caption{Geometric representation of the aspects of invariant manifold theory}
    \label{immersion 1}
\end{figure}

The notion of the convergence of the off-the-manifold dynamics to the implicit manifold and manifold attractivity can be visualized in Fig. \ref{immersion 1}. 
 The P\&I approach's concept and procedures applicable to this research are presented below.
Given a single input system
\begin{align}\label{feedbackoriginalsystem}
\begin{split}
\dt{\mathrm{x}}=\mathrm{\mathbb{F}(\mathrm{x}, \lambda)}\hspace{1.5cm}
\dt{\lambda}=\mathrm{u} 
\end{split} 
\end{align}
with $(\mathrm{x}, \lambda)\in (\mathbb{R}^{n-1},\mathbb{R})$ and without any particular structure.  The proposed P\&I approach ensures NLS stabilization by executing the four steps $(S_1-S_4)$ outlined below. \\
$(S_1)$ \textbf{Construction of the implicit manifold:}

The target dynamics  $\dt{\eta}=\beta(\eta)$  with $\mathrm{x}=\pi(\eta)=\eta x$ is defined such that the subsystem $\dt{\mathrm{x}}=\mathrm{\mathbb{F}(\mathrm{x}, \varphi (\mathrm{\mathrm{x}}))}$ for $\mathrm{C}^{\infty}$  mapping $\varphi(\mathrm{x}):\mathbb{R}^{\mathrm{n}}\rightarrow \mathbb{R}$ has a GES/GAS equilibrium at the origin by considering the relationship $\mathrm{\lambda}=\varphi (\mathrm{\mathrm{x}})$. This defines the implicit manifold $\mathbb{M}=\left \{(\mathrm{x},\lambda) \in \mathbb{R}^{\mathrm{n-1}}\times\mathbb{R} | \Psi(\mathrm{x}, \lambda):=\lambda-\varphi(\mathrm{x})=0 \right \}$, and $\pi(\eta)=\mathrm{col}(\eta, \varphi(\mathrm{\eta}))$.

{$(S_2)$}  \textbf{Tangent space structure for control systems:}
Consider an n-dimensional manifold $\mathbb{M}$ with tangent bundle $\mathbb{T}_{\mathbb{M}}$, such that all $\mathrm{p}\in {\mathbb{M}}$, $\mathbb{T}_{\mathrm{p}}{\mathbb{M}}$ has the following structure
\begin{equation}
    \mathbb{T}_{\mathrm{p}}{\mathbb{M}} = \mathbb{H}_{\mathrm{p}} \oplus  \mathbb{V}_{\mathrm{p}}: \hspace{0.3cm}  \mathbb{H}_{\mathrm{p}} \cap   \mathbb{V}_{\mathrm{p}}=0
\end{equation}
where $\mathbb{H}_{\mathrm{p}}$ is the horizontal space and $\mathbb{V}_{\mathrm{p}}$ is the vertical space.  
\begin{equation}\label{TpMm}
  \mathrm{Then} \hspace{0.9cm} \mathbb{T}_{\mathrm{p}}{\mathbb{M}}= \mathbb{H}_{\mathrm{p}} \oplus  \mathbb{V}_{\mathrm{p}}=(\dt{\mathrm{x}}, 0)\oplus (0, \dt{\lambda})=(\dt{\mathrm{x}},\dt{\lambda})
\end{equation} 
is written for given system (\ref{feedbackoriginalsystem}) at any point $\mathrm{p}\in\mathbb{M}$. With the implicit manifold $\Psi (\mathrm{x} ,\lambda )$ obtained in $S_1$, the normal vector direction is given by $\triangledown \Psi (\mathrm{x} ,\lambda )$. Thus, a PR metric $\mathrm{R}$ on space $\mathbb{T}_{\mathrm{p}}{\mathbb{M}}$ can be defined as
\begin{small}
\begin{align}\label{metric}
    \begin{split}
        \mathrm{R}&=\triangledown \Psi (\mathrm{x} ,\lambda )^{\mathrm{T}}\triangledown \Psi (\mathrm{x} ,\lambda )\\&=\begin{bmatrix}
\left ( \frac{\partial \varphi}{\partial \mathrm{x}}\right )^{\mathrm{T}}\left ( \frac{\partial \varphi}{\partial \mathrm{x}}\right ) & \left ( -\frac{\partial \varphi}{\partial \mathrm{x}}\right )^{\mathrm{T}}\\ 
-\left ( \frac{\partial \varphi}{\partial \mathrm{x}}\right ) & \mathrm{I}
\end{bmatrix}=\begin{bmatrix}
\mathrm{\mathrm{m}}_{11} & \mathrm{\mathrm{m}}_{12}\\ 
 \mathrm{\mathrm{m}}_{21}&\mathrm{\mathrm{m}}_{22}
\end{bmatrix}.
    \end{split}
\end{align}
\end{small}
which is intuitively a natural choice.
\begin{remark}
To obtain the passive output, the metric $\chi$ is replaced with semi-Riemannian metric $\mathrm{R}$ as a natural choice. 
\end{remark}
For $(\mathbb{M}, \mathrm{R})$, the splitting is visualized as follows:
\begin{small}
\begin{equation}
     (\dt{\mathrm{x}},\dt{\lambda})=\left (\dt{\mathrm{x}}, -\mathrm{m}_{22}^{-1}\mathrm{m}_{21} \dt{\mathrm{x}} \right )\oplus \left (0, \dt{\lambda}+\mathrm{m}_{22}^{-1}\mathrm{m}_{21} \dt{\mathrm{x}} \right ) =\widetilde{\mathbb{H}}_{\mathrm{p}} \oplus \widetilde{\mathbb{V}}_{\mathrm{p}}
\end{equation}
\end{small}

As $\dt{\lambda}$  is along the vertical direction, the passive output is chosen as a component $\dt{\lambda}+\mathrm{m}_{22}^{-1}\mathrm{m}_{21} \dt{\mathrm{x}}$ which  is in the same direction or parallel to $\dt{\lambda}$ \cite{nayyer2022towards}.  Roughly speaking, the idea is to bring the component $\dt{\lambda}+\mathrm{m}_{22}^{-1}\mathrm{m}_{21}$ of $\widetilde{\mathbb{H}}_{\mathrm{p}}$ to the component $\mathrm{m}_{22}^{-1}\mathrm{m}_{21} \dt{\mathrm{x}}$ of $\widetilde{\mathbb{V}}_{\mathrm{p}}$.
\begin{figure}[ht!]
    \centering
    \includegraphics[width=\linewidth]{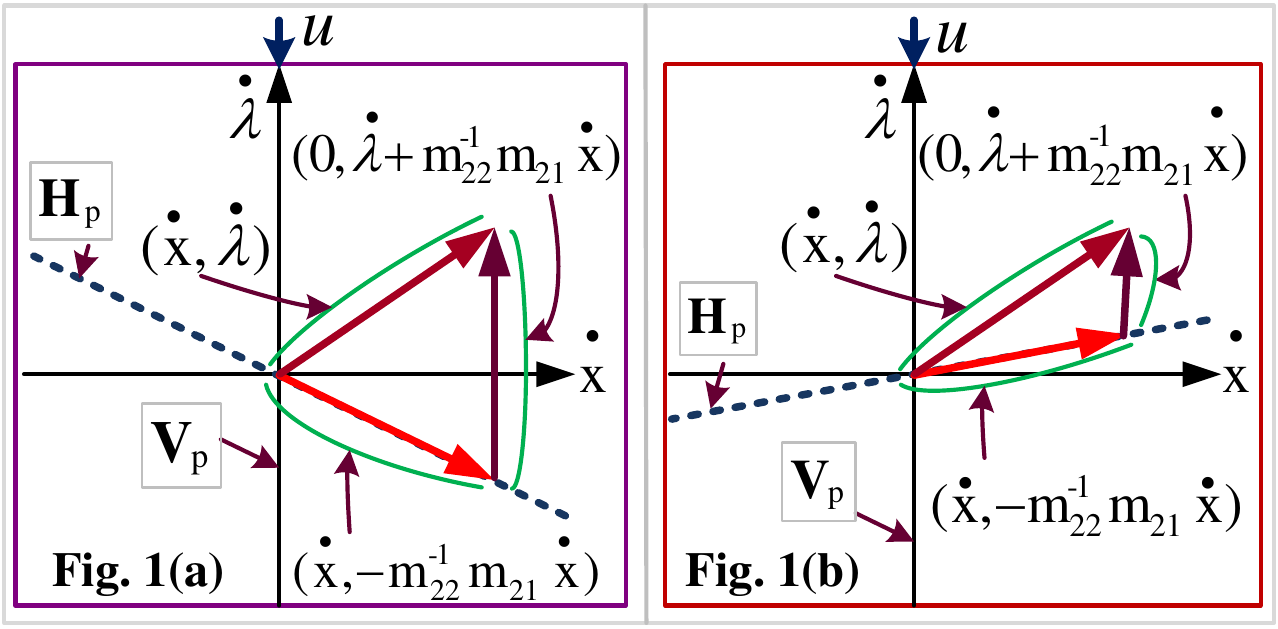}
    \caption{Geometrical interpretation: vertical vector $\mathbb{V}_{\mathrm{p}}$ is along the fiber direction and $\mathbb{H}_{\mathrm{p}} \oplus  \mathbb{V}_{\mathrm{p}}=\mathbb{T}_{\mathrm{p}}{\mathbb{M}}$. }
    \label{Geometrical interpretation}
\end{figure}
The geometrical interpretation for the splitting tangent vector is shown in Fig. \ref{Geometrical interpretation}. 

\textbf{$(S_3)$  Passive output:}
The component of $\mathrm{u}$ tangent vector along $\dt{\lambda}$  is used to define the \textit{passive output} 
\begin{equation}\label{GammaPassive}
    \mathrm{y}=\mathrm{y}_1+\mathrm{y}_2=\int_{0}^{t}\dt{\lambda}\hspace{0.08cm}\mathrm{dt}
+{\int_{0}^{t}(\mathrm{m}_{22}^{-1}\mathrm{m}_{21}\dt{\mathrm{x}})\hspace{0.08cm}\mathrm{dt}}
\end{equation}
with the help of the passivity theory. Here, $\mathrm{y}_1={\int_{0}^{t}\dt{\lambda}\hspace{0.08cm}\mathrm{dt}}$ and $\mathrm{y}_2={\int_{0}^{t}(\mathrm{m}_{22}^{-1}\mathrm{m}_{21}\dt{\mathrm{x}})\hspace{0.08cm}\mathrm{dt}}$ are defined. If $\mathrm{m}_{22}^{-1}\mathrm{m}_{21}$ is a constant then $\mathrm{y}=({\lambda}+\mathrm{m}_{22}^{-1}\mathrm{m}_{21}{\mathrm{x}})$ is defined. If $\mathrm{m}_{22}^{-1}\mathrm{m}_{21}$ is a function of $\mathrm{x}$ then, it can be written as the gradient of any function $\mathrm{q(x)}$ i.e.,  $\mathrm{m}_{22}^{-1}\mathrm{m}_{21}\mathrm{(x)}=\triangledown \mathrm{q(x)}$. Then
\begin{equation}\label{Gamma}
    \mathrm{y}=\int_{0}^{t}(\dt{\lambda}+\triangledown \mathrm{q(x)} \dt{\mathrm{x}})\mathrm{dt}=(\lambda+\mathrm{q(x)})
\end{equation}

$(S_4)$ \textbf{Storage function:}
With $\mathrm{y}$, the candidate Lyapunov function $\mathbb{S}(\mathrm{x}, \lambda)$ (i.e., storage function) is defined as
\begin{equation}\label{Storage function}
  \mathbb{S}(\mathrm{x}, \lambda)=\frac{1}{2}{\mathrm{y}}^2=\frac{1}{2}(\lambda+\mathrm{q(x)})^2.
\end{equation}
One can use the condition $ \dt{\mathbb{S}}\leq -\alpha\mathbb{S}$ along with the storage function (\ref{Storage function}) and passive output (\ref{Gamma}) to get
\begin{equation}\label{sjksdkc}
(\dt{\lambda}+\frac{\partial \mathrm{q(x)}}{\partial \mathrm{x}}\dt{\mathrm{x}})+\frac{\alpha}{2}(\lambda+\mathrm{q(x)})\leq 0
\end{equation}
The equation (\ref{sjksdkc}) is modified by substituting (\ref{feedbackoriginalsystem}) in terms of the final control law as
\begin{equation}\label{final control law}
  \mathrm{u}=-\frac{\alpha}{2}\lambda-\frac{\alpha}{2}\mathrm{q(x)}-\frac{\partial\mathrm{q(x)}}{\partial \mathrm{x}}\mathrm{f(\mathrm{x}, \lambda)}.
\end{equation}
The above-defined control law ensures the GAS equilibrium point of the system to zero/origin.
\begin{remark}
    The mathematical details \& related theorems-proofs of the P\&I approach can be found in \cite{nayyer2022towards}.
\end{remark}

\bibliography{Bibliography}

\begin{thebibliography}{10}
\providecommand{\url}[1]{#1}
\csname url@samestyle\endcsname
\providecommand{\newblock}{\relax}
\providecommand{\bibinfo}[2]{#2}
\providecommand{\BIBentrySTDinterwordspacing}{\spaceskip=0pt\relax}
\providecommand{\BIBentryALTinterwordstretchfactor}{4}
\providecommand{\BIBentryALTinterwordspacing}{\spaceskip=\fontdimen2\font plus
\BIBentryALTinterwordstretchfactor\fontdimen3\font minus
  \fontdimen4\font\relax}
\providecommand{\BIBforeignlanguage}[2]{{%
\expandafter\ifx\csname l@#1\endcsname\relax
\typeout{** WARNING: IEEEtran.bst: No hyphenation pattern has been}%
\typeout{** loaded for the language `#1'. Using the pattern for}%
\typeout{** the default language instead.}%
\else
\language=\csname l@#1\endcsname
\fi
#2}}
\providecommand{\BIBdecl}{\relax}
\BIBdecl

\bibitem{hu2009accelerated}
C.~Hu, W.~Pan, and J.~Kwok, ``Accelerated gradient methods for stochastic
  optimization and online learning,'' \emph{Advances in Neural Information
  Processing Systems}, vol.~22, 2009.

\bibitem{lan2012optimal}
G.~Lan, ``An optimal method for stochastic composite optimization,''
  \emph{Mathematical Programming}, vol. 133, no. 1-2, pp. 365--397, 2012.

\bibitem{ghadimi2016accelerated}
S.~Ghadimi and G.~Lan, ``Accelerated gradient methods for nonconvex nonlinear
  and stochastic programming,'' \emph{Mathematical Programming}, vol. 156, no.
  1-2, pp. 59--99, 2016.

\bibitem{li2015accelerated}
H.~Li and Z.~Lin, ``Accelerated proximal gradient methods for nonconvex
  programming,'' \emph{Advances in neural information processing systems},
  vol.~28, 2015.

\bibitem{beck2009fast}
A.~Beck and M.~Teboulle, ``A fast iterative shrinkage-thresholding algorithm
  for linear inverse problems,'' \emph{SIAM journal on imaging sciences},
  vol.~2, no.~1, pp. 183--202, 2009.

\bibitem{lan2011primal}
G.~Lan, Z.~Lu, and R.~D. Monteiro, ``Primal-dual first-order methods with
  iteration-complexity for cone programming,'' \emph{Mathematical Programming},
  vol. 126, no.~1, pp. 1--29, 2011.

\bibitem{polyak1964some}
B.~T. Polyak, ``Some methods of speeding up the convergence of iteration
  methods,'' \emph{Ussr computational mathematics and mathematical physics},
  vol.~4, no.~5, pp. 1--17, 1964.

\bibitem{nesterov1983method}
Y.~E. Nesterov, ``A method of solving a convex programming problem with
  convergence rate $o\left (\frac{1}{k^2} \right )$,'' in \emph{Doklady
  Akademii Nauk}, vol. 269, no.~3.\hskip 1em plus 0.5em minus 0.4em\relax
  Russian Academy of Sciences, 1983, pp. 543--547.

\bibitem{attouch2000heavy}
H.~Attouch, X.~Goudou, and P.~Redont, ``The heavy ball with friction method, i.
  the continuous dynamical system: global exploration of the local minima of a
  real-valued function by asymptotic analysis of a dissipative dynamical
  system,'' \emph{Communications in Contemporary Mathematics}, vol.~2, no.~01,
  pp. 1--34, 2000.

\bibitem{shi2021understanding}
B.~Shi, S.~S. Du, M.~I. Jordan, and W.~J. Su, ``Understanding the acceleration
  phenomenon via high-resolution differential equations,'' \emph{Mathematical
  Programming}, pp. 1--70, 2021.

\bibitem{helmke2012optimization}
U.~Helmke and J.~B. Moore, \emph{Optimization and dynamical systems}.\hskip 1em
  plus 0.5em minus 0.4em\relax Springer Science \& Business Media, 2012.

\bibitem{brockett1991dynamical}
R.~W. Brockett, ``Dynamical systems that sort lists, diagonalize matrices, and
  solve linear programming problems,'' \emph{Linear Algebra and its
  applications}, vol. 146, pp. 79--91, 1991.

\bibitem{soto2012nonlinear}
J.~Soto and J.-J.~E. Slotine, ``Nonlinear contraction tools for constrained
  optimization,'' \emph{arXiv preprint arXiv:1206.1838}, 2012.

\bibitem{lessard2016analysis}
L.~Lessard, B.~Recht, and A.~Packard, ``Analysis and design of optimization
  algorithms via integral quadratic constraints,'' \emph{SIAM Journal on
  Optimization}, vol.~26, no.~1, pp. 57--95, 2016.

\bibitem{wibisono2016variational}
A.~Wibisono, A.~C. Wilson, and M.~I. Jordan, ``A variational perspective on
  accelerated methods in optimization,'' \emph{proceedings of the National
  Academy of Sciences}, vol. 113, no.~47, pp. E7351--E7358, 2016.

\bibitem{9971723}
S.~Shadab~Nayyer, S.~R. Wagh, and N.~M. Singh, ``Passivity and immersion (p\&i)
  approach for constructive stabilization and control of nonlinear systems,''
  \emph{IEEE Control Systems Letters}, vol.~7, pp. 817--822, 2023.

\bibitem{nayyer2022passivity_gp}
S.~S. Nayyer, J.~Hozefa, R.~Gunjal, S.~Bhil, S.~Wagh, and N.~Singh, ``Passivity
  and immersion (p\&i) approach with gaussian process for stabilization and
  control of nonlinear systems,'' \emph{IEEE Access}, vol.~10, pp.
  132\,621--132\,634, 2022.

\bibitem{nayyer2022passivity}
S.~S. Nayyer, G.~Revati, S.~Wagh, and N.~Singh, ``Passivity and immersion
  based-modified gradient estimator: A control perspective in parameter
  estimation,'' \emph{arXiv preprint arXiv:2211.10674}, 2022.

\bibitem{wensing2020beyond}
P.~M. Wensing and J.-J. Slotine, ``Beyond convexity—contraction and global
  convergence of gradient descent,'' \emph{Plos one}, vol.~15, no.~8, p.
  e0236661, 2020.

\bibitem{astolfi2008nonlinear}
A.~Astolfi, D.~Karagiannis, and R.~Ortega, \emph{Nonlinear and adaptive control
  with applications}.\hskip 1em plus 0.5em minus 0.4em\relax Springer, 2008,
  vol. 187.

\bibitem{nayyer2022towards}
S.~S. Nayyer, S.~R. Wagh, and N.~M. Singh, ``Towards a constructive framework
  for stabilization and control of nonlinear systems: Passivity and immersion
  (p$\backslash$\&i) approach,'' \emph{arXiv preprint arXiv:2208.10539}, 2022.

\bibitem{alvarez2002second}
F.~Alvarez, H.~Attouch, J.~Bolte, and P.~Redont, ``A second-order gradient-like
  dissipative dynamical system with hessian-driven damping.: Application to
  optimization and mechanics,'' \emph{Journal de math{\'e}matiques pures et
  appliqu{\'e}es}, vol.~81, no.~8, pp. 747--779, 2002.

\bibitem{talasila2006normally}
V.~Talasila, J.~Clemente-Gallardo, and A.~Astolfi, ``Normally hyperbolic
  controlled-invariant manifolds,'' in \emph{17th International Symposium on
  Mathematical Theory of Networks and Systems}, 2006.

\bibitem{fenichel1971persistence}
N.~Fenichel and J.~Moser, ``Persistence and smoothness of invariant manifolds
  for flows,'' \emph{Indiana University Mathematics Journal}, vol.~21, no.~3,
  pp. 193--226, 1971.

\bibitem{hirsch1977invariant}
M.~Hirsch, C.~Pugh, and M.~Shub, ``Invariant manifolds (lecture notes in
  mathematics, 583),'' 1977.

\bibitem{capinski2020persistence}
M.~J. Capi{\'n}ski and H.~Kubica, ``Persistence of normally hyperbolic
  invariant manifolds in the absence of rate conditions,'' \emph{Nonlinearity},
  vol.~33, no.~9, p. 4967, 2020.

\bibitem{wang2015immersion}
L.~Wang, F.~Forni, R.~Ortega, and H.~Su, ``Immersion and invariance
  stabilization of nonlinear systems: A horizontal contraction approach,'' in
  \emph{2015 54th IEEE Conference on Decision and Control (CDC)}.\hskip 1em
  plus 0.5em minus 0.4em\relax IEEE, 2015, pp. 3093--3097.

\bibitem{nesterov2003introductory}
Y.~Nesterov, \emph{Introductory lectures on convex optimization: A basic
  course}.\hskip 1em plus 0.5em minus 0.4em\relax Springer Science \& Business
  Media, 2003, vol.~87.

\bibitem{sun2020high}
B.~Sun, J.~George, and S.~Kia, ``High-resolution modeling of the fastest
  first-order optimization method for strongly convex functions,'' in
  \emph{2020 59th IEEE Conference on Decision and Control (CDC)}.\hskip 1em
  plus 0.5em minus 0.4em\relax IEEE, 2020, pp. 4237--4242.

\end{thebibliography}
\bibliographystyle{IEEEtran}


\end{document}